\documentclass[journal]{IEEEtran}
\usepackage{cite}
\ifCLASSINFOpdf
   \usepackage[pdftex]{graphicx}
   \graphicspath{{../pdf/}{../jpeg/}}
   \DeclareGraphicsExtensions{.pdf,.jpeg,.png}
\else
  \usepackage[dvips]{graphicx}
  \graphicspath{{../eps/}}
  \DeclareGraphicsExtensions{.eps}
\fi
\usepackage[outdir=./sim15Results/]{epstopdf}
\usepackage{amssymb} 
\usepackage{amsmath} 
\usepackage{bbold}
\usepackage{array}
\ifCLASSOPTIONcompsoc
 \usepackage[caption=false,font=normalsize,labelfont=sf,textfont=sf]{subfig}
\else
  \usepackage[caption=false,font=footnotesize]{subfig}
\fi
\begin{document}

\title{Resilient Mobile Multi-Target Surveillance Using Multi-Hop Autonomous UAV Networks for Extended Lifetime}

\author{Abdulsamet~Dağaşan,
        and~Ezhan~Karaşan
\thanks{A. Dağaşan and E. Karaşan with the Department
of Electrical and Electronics Engineering, Bilkent University, Ankara
06800, Turkey.}
}

%



\maketitle

\begin{abstract}
Cooperative utilization of Unmanned Aerial Vehicles (UAVs) in public and military surveillance applications has attracted significant attention in recent years. Most UAVs are equipped with sensors that have bounded coverage and wireless communication equipment with limited range. Such limitations pose challenging problems to monitor mobile targets. This paper examines fulfilling surveillance objectives to achieve better coverage while building a resilient network between UAVs with an extended lifetime. The multiple target tracking problem is studied by including a relay UAV within the fleet whose trajectory is autonomously calculated in order to achieve a reliable connected network among all UAVs. Optimization problems are formulated for single-hop and multi-hop communications among UAVs. Three heuristic algorithms are proposed for multi-hop communications and their performances are evaluated. A hybrid algorithm, which dynamically switches between single-hop and multi-hop communications is also proposed. The effect of the time horizon considered in the optimization problem is studied. Performance evaluation results show that the trajectories generated for the relay UAV by the hybrid algorithm can achieve network lifetimes that are within 5\% of the maximum possible network lifetime which can be obtained if the entire trajectories of all targets were known a priori. \end{abstract}

\begin{IEEEkeywords}
UAVs, multi-target surveillance, resilient multi-hop network topology, network lifetime.
\end{IEEEkeywords}

%
\IEEEpeerreviewmaketitle

\section{Introduction}
%
%
%
%
\IEEEPARstart{S}{urveillance} systems have recently received a significant attention due to the rapid increase in security and safety threats. Some of the surveillance applications include search and rescue operations, monitoring an environment and tracking mobile targets. Among the available surveillance methods, use of Unmanned Aerial Vehicles (UAV) for these tasks has been rather widespread for mainly two reasons \cite{ref1}. Firstly, UAVs can operate where it might be too dangerous for humans to fulfill surveillance duties. Secondly, UAVs can also be operated autonomously without human resource allocation. 
The rapid advancement in UAV technology has propelled the significance of utilizing cooperative UAVs for tasks related to surveillance \cite{ref2}, tracking \cite{ref3}, and mobile target search \cite{ref4}. This significance is evident in various applications, particularly in search and rescue operations, where cooperative UAVs are employed to find and observe targets, in addition to coordinating their flight paths \cite{ref5, ref6, ref7, ref8, ref9}.

The objective of target tracking with cooperative UAVs is to monitor a specific environment and acquire information on the mobile targets. This information is analyzed locally and shared either with a central node or among UAVs. Since the communication range of UAV is limited, connectivity of the UAVs becomes an issue as the mobile targets move away from each other. Therefore, building a reliable and connected network between UAVs for achieving an extended lifetime is critical. 

Effective coordination is fundamental to the success of cooperative UAVs in target tracking. UAV coordination involves methods such as centralization, where a central station manages task distribution and path planning, decentralization, where leader UAVs collaborate on tasks and trajectories, and distributed coordination, allowing each UAV independent task and trajectory determination \cite{tortonesi2012multiple,ref10, ref11, ref12, 9926144}. Additionally, a relay UAV can be used for a succesful coordination amongst UAVs \cite{9915455}. The trajectory planning for relay UAVs is essential to optimize specific objectives, often involving sophisticated techniques like Gaussian process models, nonlinear model predictive control, and genetic algorithms \cite{ref13, ref14, ref15, giagkos2021uav, ref16}. Various studies explore coordination strategies, such as optimizing network coverage and selecting relay UAVs. Resource optimization and multi-channel, multi-radio models are also under consideration. Particularly in disaster surveillance scenarios, UAV-enabled communication systems, featuring multi-hop and innovative routing algorithms, are developed to enhance overall system performance, specifically for search and rescue operations \cite{ref16, liu2019task}.

This paper proposes trajectory planning algorithms that ensure a connected and robust network topology among distributed cooperative target tracking UAVs with the assistance of a relay UAV. A connected network topology indicates that every UAV in the network has a way of communicating with each other. This study considers both single-hop and multi-hop communication techniques to perform path planning for the relay UAV. In single-hop communication, the relay UAV is directly connected to other drones, creating a star topology. On the other hand, in multi-hop communication, the relay UAV can communicate with each UAV through other drones, resulting in an arbitrary tree topology. The objective is to maximize the duration for which a connected network topology can be sustained among UAVs.

A comparative analysis of single-hop and multi-hop communication techniques for UAVs is presented. A probabilistic coverage model which takes drone position errors into account is proposed. The results of the single-hop communication model are compared with the \textit{Center of Mass} approach as a baseline where the relay UAV is positioned at the center of mass of the tracker UAVs. It is demonstrated that the proposed method 16\% improvement when time horizon is used. Furthermore, three algorithms for multi-hop communication are proposed. The network topologies generated by the three proposed algorithms remain connected for over 90\% of the maximum possible network lifetime whereas the most effective algorithm among the three proposed algorithms achieve network lifetimes that are within less than 5\% of the maximum lifetime that can be obtained if all target locations were known a priori. Overall, the potential of using UAVs for both single-hop and multi-hop communication scenarios is demonstrated, and insights into effective autonomous planning strategies for optimizing the performance of such systems are provided.

The organization of the rest of the article is as follows: Section II provides an overview of the existing research in the field. Section III introduces the system model, position estimator and path planner of relay drone. In Section IV, the problem formulation of trajectory planning is explained. Section V illustrates the numerical results for the estimator and path planner, respectively. The theoretical limits on network lifetime are discussed in Section VI and the performances of proposed algorithms are compared with these limits. Lastly, Section VII concludes the paper.

\section{Related Work}
Coordination amongst the UAVs is critical to achieve successful cooperative target tracking. Centralized coordination means that UAVs transmit the information to a central station. Task distribution and path planning for every UAV is calculated at the central station and sent back to the on-duty UAVs \cite{tortonesi2012multiple}. Decentralized coordination means there are multiple UAVs that act as central nodes. The information from the other UAVs are received by the leader UAVs, then they cooperatively assign tasks and plan trajectories for the on-duty UAVs \cite{ref10}. The UAVs can also coordinate in a distributed manner. Distributed coordination means each UAV decides its task and plans its trajectory independently \cite{ref11}, \cite{ref12}, \cite{9926144}. 

Relay UAVs which designed to act as an intermediary bridge between other UAVs or ground-based communication devices are used to extend the communication range between devices that may be out of direct range of each other. Regardless of the coordination method, a relay UAV can be used to achieve reliable target tracking \cite{9915455}. For the trajectory planning of relay UAV, the task is to plan the path which optimizes a given objective. \cite{ref13} investigates path planning for a relay UAV in urban areas for Airborne to Ground (A2G) communications where the ground nodes are stationary. A Gaussian process model is formed and it is solved with a Nonlinear Model Predictive Control (NMPC)-based trajectory planner. In \cite{ref14}, the same problem is studied with mobile ground nodes. NMPC-based trajectory planner is also used in this study, however, a discrete genetic algorithm is used to find the optimal control input. To find the optimal path,  NMPC is combined with a finite time horizon controlled system. Time horizon determines the time window that the relay UAV movements will be optimized. \cite{ref15} uses the same NMPC-based planner, but to calculate the objective function, it uses the concept of Minimum Spanning Tree (MST). MST is used to choose the communication link with the highest expectation of a successful transmission. Group of autonomous UAVs to maximize the network coverage of mobile ground targets are examined in \cite{giagkos2021uav} using the game theory, and the performance is compared with the genetic algorithm.

Selection of relay UAVs from a set of cooperative UAVs is suggested in \cite{ref16}. A matching market-based optimization for UAV communication models are presented for different coordination among cooperative UAVs in order to choose relay UAVs. \cite{ref16} analyzes UAV networks from the perspective of resource optimization. The idea is further extended in \cite{liu2019task} and multi-channel, multi-radio competition model is modeled with different objective functions. 

Communication coverage amongst UAVs with communication range restrictions and requirement of multi-hop communication between UAVs and a base station are studied in \cite{sharvari2023connectivity}. The study focuses on the routing problem in a communication system enabled by UAVs for disaster surveillance. The system includes multiple UAVs that communicate with each other and a remote Terrestrial Base Station (TBS), with the aim of exploring the disaster-affected areas, collecting information on the event's impact, and relaying the data to aid search and rescue operations. A novel algorithm called the Multihop Opportunistic 3D Routing (MO3DR) is proposed, which addresses coverage and collision constraints without the need for trajectory planning, to facilitate efficient data transmission in the system. The aim of the study is to improve the performance of the UAV-enabled communication in disaster scenarios. 

The proposed algorithm in \cite{YANMAZ2022102800} is a modular approach for positioning relays and planning trajectories for UAV missions. It ensures that the UAV mission team maintains connectivity with minimal relays and feasible paths. The study presents different strategies for relay positioning based on the objective and compares them to ideal and Voronoi-based benchmark schemes. The concept is further developed in \cite{9430858} by introducing dynamic relay selection from a group of available UAVs.

To the best of our knowledge, previous research has not explored the concept of ensuring reliable communication with distributed cooordination through the use of a relay UAV that autonomously plans its trajectory based on single-hop and multi-hop network topology scenarios, while also utilizing a probabilistic coverage model which takes the estimation errors in UAV positions into account.

\section{System Model}
To achieve connectivity between $N$ drones that are tracking targets, a relay UAV is utilized to enhance the connectivity among the drones. The relay UAV autonomously decides its path to communicate with other drones, which comprises three components: sensors, estimator, and path planner. The sensor and estimator components primarily aim to develop an optimization-based path planner.

Each tracker drone is equipped with Inertial Measurement Units (IMU) sensors and a Global Navigation Satellite System (GNSS), while the relay drone has a GNSS receiver. Using the data from the sensors, the relay drone estimates the positions of the drones by considering the position and velocity information from the sensors to estimate the position and path of the other drones. To estimate the position and velocity over time, the Kalman Filter is utilized. 

\subsection{Position Estimator}
Let $\mathbf{p} = [p_x,p_y,p_z]^T$ denote the position vector and $\mathbf{v} = [v_x,v_y,v_z]^T$ denote the velocity vector as the state vector. The acceleration vector, denoted by $\mathbf{a} = [a_x,a_y,a_z]^T$, is utilized as the control vector. The state vector can be written as:
\begin{equation}
\label{eq:1}
\mathbf{x} = 
\begin{bmatrix}
    \mathbf{p} \\
    \mathbf{v} 
\end{bmatrix}.
\end{equation}
The state prediction stage can be written as:
\begin{align}
\label{eq:2}
\mathbf{x_k} = \begin{bmatrix}
\mathbf{p_k} \\
\mathbf{v_k} 
\end{bmatrix}
=
\begin{bmatrix}
\mathbf{\mathbf{p_{k-1}}} + \mathbf{\mathbf{v_{k-1}}}\Delta t + \frac{1}{2}\mathbf{\widetilde{a}_{k-1}}\Delta t^2 \\
\mathbf{v_{k-1}} + \mathbf{\widetilde{a}_{k-1}}\Delta t
\end{bmatrix}.
\end{align}

The state prediction equation can be rearranged by employing identity and zero matrices to form the following system:
\begin{align}
\label{eq:3}
\mathbf{x_k} = \begin{bmatrix}
\mathbf{I} & \mathbf{I} \Delta t\\
\mathbf{0} & \mathbf{I}  
\end{bmatrix}
\mathbf{x_{k-1}} + 
\begin{bmatrix}
\frac{1}{2} \mathbf{I} \Delta t^2 \\
\mathbf{I}\Delta t
\end{bmatrix}
\mathbf{\widetilde{a}_{k-1}}.
\end{align}

The process noise on the acceleration vector can be written with its zero mean Gaussian noise vector $\mathbf{e}$. 
\begin{align}
\label{eq:4}
\mathbf{a_{k-1}} = \mathbf{\widetilde{a}_{k-1}} + \mathbf{e_{k-1}}, \\
\mathbf{e_{k-1}} \sim N(0,\mathbf{I} \mathbf{\sigma _e^2}).
\end{align}

The process noise covariance matrix can be calculated using this noise:
\begin{equation}
\label{eq:6}
\mathbf{Q} = 
\begin{bmatrix}
\frac{1}{4} \mathbf{I} \Delta t^4 & \mathbf{0} \\
\mathbf{0} & \mathbf{I}\Delta t^2
\end{bmatrix}
\mathbf{\sigma _e^2}.    
\end{equation}

The system model is given by:
\begin{equation}
\label{eq:7}
    \mathbf{x_k} =  \mathbf{F} \mathbf{x_{k-1}} + \mathbf{B} \mathbf{a_{k-1}} + \mathbf{w_{k-1}},
\end{equation}

where,
\begin{align}
\label{eq:8}
\mathbf{F} &= \begin{bmatrix}
\mathbf{I} & \mathbf{I} \Delta t\\
\mathbf{0} & \mathbf{I}  
\end{bmatrix},
\\
B &= \begin{bmatrix}
\frac{1}{2} \mathbf{I} \Delta t^2 \\
\mathbf{I}\Delta t
\end{bmatrix}^T,
\\
    \mathbf{w_{k-1}} &\sim N(0,\mathbf{Q}).
\end{align}

The following equation can be written for the measurement vector $\mathbf{z_k}$:
\begin{equation}
\label{eq:11}
\mathbf{z_k} = \begin{bmatrix}
\mathbf{p_k} \\
\mathbf{v_k} \\
\end{bmatrix} + \mathbf{v_k},
\end{equation}
where $\mathbf{v_k} \sim N(0,\mathbf{R})$. Using $\mathbf{z_k} = \mathbf{H} \mathbf{x_k} + \mathbf{v_k}$ and \eqref{eq:11}, $\mathbf{H} = \mathbf{I_{6x6}}$ \cite{ref33} is obtained. 

The tracker drone positions can be estimated using the Kalman Filter. The estimates, along with a probability parameter, will be utilized to compute the cost function in the optimization problem.

\subsection{Path Planner}
In order to create a connected network between the drones, a dynamic optimization problem must be solved by the relay drone to determine its own path autonomously. The aim is to create a connected network within a designated time frame denoted as $[0, T]$. However, solving the optimization problem at the beginning is not feasible due to the high computational complexity and time constraints, particularly for larger values of T. Moreover, it is essential to consider environmental changes, target paths, and drone paths. Therefore, a receding horizon approach to solve the optimization problem is proposed.

To estimate the drone positions within a certain time frame $[t, t + t_k]$, where $t_k = K \times t_s$, an estimator with a fixed step size of $t_s$ is utilized. A path planning strategy is then created for a specific duration of $[t, t + t_f]$, where the optimization horizon window's duration, denoted as $t_f = F \times t_s$, is maximized to ensure that the drones remain connected for as long as possible. The selection of $F$ and $K$ is such that they are positive integers ($F,K \in \mathbf{Z^{0+}}$) and that $F$ is a factor of $K$ ($F\mid K$). If $F = 0$, the optimization problem is referred to as the \textit{no horizon problem}, as the objective function is solved only to obtain the current time solution.  

A constant-altitude kinematic model is used for the relay drone, considering its physical constraints. The state parameter of the drone is denoted by $\mathbf{x(t)} = [\mathbf{p(t)},\mathbf{v(t)}]$, and the control input parameters are represented by $\mathbf{u(t)} = [V_D(t),\theta(t)]$, which correspond to the velocity and angle of the drone, respectively. The velocity of the drone is limited to a minimum value of $V_{min}$ and a maximum value of $V_{max}$, while the change in the drone's angle for each step size is constrained to the interval $[-\theta_{max}, \theta_{max}]$. The dynamics of the drone are expressed as follows:
\begin{equation}
\label{eq:13}
V_x(t) = V_D(t) \cos{(\theta(t))}, \\
\end{equation}
\begin{equation}
\label{eq:14}
V_y(t) = V_D(t) \sin{(\theta(t))}, \\
\end{equation}
\begin{equation}
\label{eq:15}
|\theta(t)-\theta(t-1)| \leq \theta_{max},\\
\end{equation}
\begin{equation}
\label{eq:16}
V_{min} \leq V_D(t) \leq V_{max}.
\end{equation}

The aim of the optimization problem is to find the optimal input vector $\mathbf{u^{*}(t)} \in U$ which maximizes the probability that the network topology is connected at time t. $U$ represents the set of possible inputs defined in \eqref{eq:15} and \eqref{eq:16}. For arbitrary time $t_0$, given the variable $\mathbf{u^{*}(t)}$ for the interval $t \in [t_0, t_0+t_f]$ the optimal state parameter $\mathbf{x^{*}(t)}$ can be found. 

The dynamic optimization problem can be formulated similar to a Bolza-type Optimal Control Problem (OCP) \cite{ref34}. The Lagrange term can be defined as:
\begin{equation}
\label{eq:17}
\Phi_L(t,u) = \rho (t,\mathbf{x(t)},\mathbf{u(t)}), 
\end{equation}
where $\rho (t,\mathbf{x(t)},\mathbf{u(t)})$ denotes the probability of having a connected network topology with state $\mathbf{x_t}$ at time $t$. The Mayer term can be written as:
\begin{equation}
\label{eq:18}
\Phi_M(t) = d(\mathbf{x(t_0 +t_f)},\mathbf{x_t(t_0 +t_f)}),
\end{equation}
where $d(\mathbf{x(t_0 + t_f)},\mathbf{x_t(t_0 + t_f)})$ is a function of distance between drones for at the end of the optimization horizon for $t \in [t_0, t_0 +t_f]$. The resulting optimization problem is given by:
\begin{equation}
\begin{array}{ll@{}ll}
\underset{u(t)}{\text{max}}  & \displaystyle \rho (t,\mathbf{x(t)},\mathbf{u(t)}) + d(\mathbf{x(t_0 + t_f)},\mathbf{x_t(t_0 + t_f)}) \\
\text{s.t.}& (\ref{eq:13}), (\ref{eq:14}), (\ref{eq:15}), (\ref{eq:16})\\
& t \in [t_0, t_0 + t_f]
\end{array}
\label{eq:19}
\end{equation}

The objective of the task is to select the optimal values of $V_D(t)$ and $\theta(t)$ for a given interval $t \in [t_0, t_0 + t_f]$ that maximize the objective function. In order to assess the effect of the optimization horizon, the optimization problem is also solved for $F = 0$, which corresponds to $t_f = 0$. The objective function $\Phi_L(t,u) + \Phi_M(t)$ varies depending on the type of the network: single-hop or multi-hop communication. 

\section{Problem Formulation}
\subsection{Single-Hop Communication}
In this section, tracker drones and the relay drone should communicate with a single-hop communication. The velocity and position vector of $i^{th}$ drone is expressed as $\mathbf{v_i(t)}$ and $\mathbf{p_i(t)}$ at time $t$. The position of the relay drone is $\mathbf{p(t)}$ at time $t$. 

All the tracker drones should be connected to the relay drone directly. To check the connectivity between relay drone and the tracker drones, distance parameter is defined as:
\begin{equation}
\label{eq:20}
d_i(t) = || \mathbf{p(t)} - \mathbf{p_i(t)}||_2,
\end{equation}
where $d_{i}(t)$ is the distance between $i^{th}$ tracker drone and relay drone at time $t$.

\begin{figure}[t]
    \centering
    \includegraphics[width=0.6\linewidth]{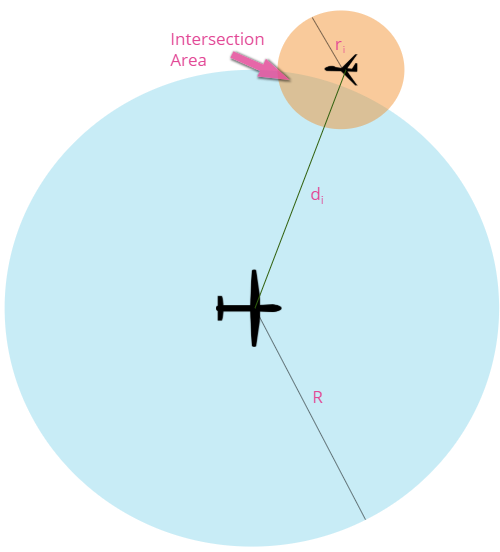}
    \caption{Blue circle is the communication area of the relay drone and the orange circular region is the estimated drone position area.}
    \label{fig:Intersection}
\end{figure}

In order to model the inaccuracy in estimating the velocity and position of drones, a probabilistic model is used as opposed to a deterministic model. A probable convex region is assigned to the estimated drone to guarantee the maximization of the connectivity between drones. In this case, a circular region is assigned with radius $r_i(t) = k \cdot \mathbf{v_i(t)}$ and center $\mathbf{p_i(t)}$ to a tracker drone. Assuming a uniformly distributed position error, the intersection area between the circular drone position uncertainty region and the area corresponding to the communication range is calculated to get the connectivity likelihood function $I(R,r,l)$. The visualization for an example scenario is shown in Figure \ref{fig:Intersection}.

The intersection area $I(R,r,d)$ between two circles can be found as:
\begin{multline}
\label{eq:21}
I(R,r,d) = \\
R^2 \arccos{\bigg(\frac{R^2 - r^2 + d^2}{2 d R}\bigg)} + r^2 \arccos{\bigg(\frac{r^2 - R^2 + d^2}{2 d r}\bigg)} \\
-\frac{1}{2}\bigg((R + r + d)(R - r + d)(-R + r + d)(R + r - d)\bigg)^{\!1/2},
\end{multline}
where $R$ and $r$ are the radii of two circles and $d$ is the distance between the centers of two circles. 

The probability of connectivity for drone i, denoted as $\rho_i$, is expressed as the ratio of the intersection area between the circular regions of the tracker drone and the relay UAV, to the total area of the circular region of the tracker drone. 
\begin{equation}
\label{eq:22}
\rho_i (t,\mathbf{x(t)},\mathbf{u(t)}) =
 \begin{cases}
\frac{I(R,r_{i}(t),d_{i}(t))}{\pi r_i(t)^2},& \text{if } d_{i}(t) \leq R + r_i(t)\\
 0,              & \text{otherwise}.\\
 \end{cases}
\end{equation}

By using the independence assumption for connectivity of individual drones and position errors between different UAVs, the overall probabilistic connectivity function $\rho$ can be written in the product form. The probability that all tracker drones are connected to the relay UAV can be expressed as follows:
\begin{equation}
\label{eq:23}
\rho (t,\mathbf{x(t)},\mathbf{u(t)}) = \prod_{i = 1}^{N-1} \rho_i (t,\mathbf{x(t)},\mathbf{u(t)}),  
\end{equation}
where $N$ is the total number of drones, including the relay drone. 

In single-hop communication, it is necessary for the relay drone to establish a direct connection with the tracker drones. The proposed approach involves closely monitoring each tracker drone to accommodate any changes in their movements. For this purpose,  $d(\mathbf{x(t)},\mathbf{x_t(t)})$, is utilized which aims to minimize the maximum distance between the relay drone and the tracker drones. In order to avoid this component dominating the objective function, it is multiplied by a small value of $\epsilon$, $0 < \epsilon << 1$. The function can be expressed as follows:
\begin{equation}
\label{eq:24}
    d(\mathbf{x(t)},\mathbf{x_t(t)}) = - \epsilon \cdot \max \{d_i(t)\}.
\end{equation}

The resulting optimization problem is expressed as: 

\begin{equation}
\label{eq:25}
\begin{array}{ll@{}ll}
\underset{u(t)}{\text{max}}  & \displaystyle \prod_{i}^{N-1} \rho_i (t,\mathbf{x(t)},\mathbf{u(t)}) - \epsilon \cdot \max \{d_i(t_0 + t_f)\}\\
\text{s.t.}& \eqref{eq:13}, \eqref{eq:14}, \eqref{eq:15}, \eqref{eq:16}\\
& t \in [t_0, t_0 + t_f]
\end{array}
\end{equation}

While the primary objective is to maximize the probability that the relay is connected to all tracker UAVs, the secondary objective is used to make sure that among all solutions that optimize the primary objective, the one which minimizes the maximum distance between the relay and tracker drones is chosen. The path planner component solves this optimization problem with the given position estimates of the other drones in a chosen horizon window and generates the trajectory. 

\subsection{Multi-Hop Communication}
Within this section, multi-hop communication is utilized by the tracker drones and the relay drone, meaning that direct connection between them is not required. Instead, the goal is to establish a connected network topology among all N drones. To solve the dynamic optimization problem, the objective function defined in \eqref{eq:19} is employed. The Lagrange term $\rho(t, \mathbf{x(t)}, \mathbf{u(t)})$ and the Mayer term $d(\mathbf{x(t_0 +t_f)},\mathbf{x_t(t_0 +t_f)})$ of the objective function $\Phi_L(t,u) + \Phi_M(t)$ will differ from those used in the context of single-hop communication.

Multi-hop communication is established by creating a network graph consisting of $N$ drones, including the relay drone, where each drone is represented as a node in the graph. For all possible $\binom{N}{2}$ edges, a state matrix $\mathbf{S}$ is defined to calculate the probabilistic connectivity function $\rho$. The state matrix contains all possible combinations of edges that can be connected. Consequently, the size of the state matrix $\mathbf{S}$ is $\binom{N}{2} \times 2^{\binom{N}{2}}$. To calculate the probability of a connected network, the connectivity probability for each edge is computed, and then multiplied with the state matrix to obtain the value of the probabilistic connectivity function $\rho$.

The connectivity probability for each edge is determined by calculating the average intersection area between each drone, utilizing the same model depicted in Figure \ref{fig:Intersection}. For edges connecting the relay drone to other drones, the probability of connection derived in \eqref{eq:22} is used. However, for network edges connecting two tracker drones, a different probabilistic model is employed due to position inaccuracy being valid for both drones. Thus, an average intersection area between tracker drones $i$ and $j$ is computed to determine the probability of establishing a link between UAVs $i$ and $j$. 

The Probability Density Function (PDF) of distance between two tracker drones, $f(x;d,r)$, can be geometrically formulated as: 
\begin{equation}
\label{eq:26}
f(x;d,r) = \frac{2x}{\pi r^2}\arccos{\bigg(\frac{x^2-r^2+d^2}{2d x}\bigg)}.
\end{equation}

The average intersection area can be found as follows:

\begin{equation}
\label{eq:27}
I_{ij}(t) = \int_{d_{ij}(t)-r_i(t)}^{d_{ij}(t)+r_i(t)} I(R_c,r_{i}(t),x)f(x;d_{ij}(t),r_j(t)) dx
\end{equation}
where $I(R_c,r_{i}(t),x)$ is the intersection area between drone $j$ with communication radius $R_{c}$ and drone $i$ with location region with radius $r_i(t)$ at time $t$. It is assumed that all tracker drones possess the same communication radius $R_c$. The average intersection area is calculated by integrating over the distance $d_{ij}(t) - r_{j}(t) \leq x \leq d_{ij}(t) + r_{j}(t)$ between two drones.

By applying \eqref{eq:27} to each pair of vertices in the matrix, a vector of size 1$ \times \binom{N}{2}$ is derived. This vector is multiplied with the state matrix $\boldsymbol{S}$ to obtain a vector of size 1$\times2^{\binom{N}{2}}$ which contains the probabilities for all possible network configurations. The transpose of this vector is multiplied with the one vector $\mathbb{1}_{1 \times 2^{\binom{N}{2}}}$ to determine the probabilistic connectivity function.
\begin{equation}
\label{eq:28}
\rho (t,\mathbf{x(t)},\mathbf{u(t)}) = \mathbb{1}(\mathbf{I(t)S})^T.
\end{equation}

The relay drone is responsible for establishing a connected network in multi-hop communication. Three algorithms are proposed to enable the relay drone to strategically position itself and sustain connectivity for longer periods while drone positions change over time.
\subsubsection{Nearest Point Algorithm}
The first algorithm involves dividing the drones into two sets, denoted as $S_1$ and $S_2$. The drone furthest away from the other drones is placed in $S_1$, while the rest of the drones are in $S_2$. The distances of all drones in both sets to the relay drone are measured, and the closest ones are selected. Let drone $k \in S_1$ and drone $l \in S_2$ be the closest to the relay drone, with distances $d_k(t)$ and $d_l(t)$, respectively. In the next iteration, the control input that minimizes $\max {d_k(t), d_l(t)}$ is chosen. To prevent this part from dominating the objective function, it is multiplied by a small $\epsilon > 0$ value. This can be used as the Mayer part of the optimal control problem as follows:
\begin{equation}
\label{eq:29}
d(\mathbf{x(t)},\mathbf{x_t(t)}) = -\epsilon \cdot \max{\{d_{k}(t),d_{l}(t)\}}
\end{equation}

The objective function of the resulting optimization problem can be written as:
\begin{equation}
\label{eq:30}
\begin{array}{ll@{}ll}
\text{max}  & \displaystyle \mathbb{1}(\mathbf{I(t)S})^T
- \epsilon \cdot \max{\{d_{k}(t_0 +t_f),d_{l}(t_0 +t_f)\}}\\
\text{s.t.}& (\ref{eq:13}), (\ref{eq:14}), (\ref{eq:15}), (\ref{eq:16})\\
& t \in [t_0, t_0 + t_f]
\end{array}
\end{equation}

\subsubsection{Midpoint Algorithm}
In this algorithm, the centroid of the positions of drones in $S_2$ is first calculated and denoted by $\mathbf{p_{mean}(t)}$. Next, the distance between the relay drone and $\mathbf{p_{mean}(t)}$ is calculated and denoted as $d_{mean}(t)$.
\begin{equation}
\label{eq:31}
d_{mean}(t) = || p(t) - \mathbf{p_{mean}(t)}||_2.
\end{equation}

The control input is selected to minimize $\max {d_{mean}(t), d_k(t)}$ in this algorithm. To prevent this part from dominating the objective function, it is multiplied by a small $\epsilon > 0$ value. This can be utilized as the Mayer part of the OCP as follows:
\begin{equation}
\label{eq:33}
d(\mathbf{x(t)},\mathbf{x_t(t)}) = -\epsilon \cdot \max{\{d_{k}(t),d_{mean}(t)\}}
\end{equation}

The objective function of the optimization problem can be expressed as:
\begin{equation}
\label{eq:34}
\begin{array}{ll@{}ll}
\text{max}  & \displaystyle \mathbb{1}(\mathbf{I(t)S})^T- \epsilon \cdot \max{\{d_{k}(t_0 +t_f),d_{mean}(t_0 +t_f)\}}\\
\text{s.t.}& (\ref{eq:13}), (\ref{eq:14}), (\ref{eq:15}), (\ref{eq:16})\\
& t \in [t_0, t_0 + t_f]
\end{array}
\end{equation}

\subsubsection{Hybrid Algorithm}
A hybrid approach is employed in the last algorithm, which dynamically switches between single-hop and multi-hop objectives. The algorithm prioritizes single-hop communication and solves the optimization problem introduced in \eqref{eq:25} if single-hop communication is possible. If the obtained solution does not satisfy the following conditions:
\begin{equation}
    \underset{u(t)}{\text{max}}  \displaystyle \prod_{i}^{N-1} \rho_i (t,\mathbf{x^{*}(t)},\mathbf{u^{*}(t)}) \geq 0,
\end{equation}
single-hop connectivity is not achieved, hence the algorithm switches to multi-hop communication optimization problem. The optimization problem in \eqref{eq:34} is then solved to obtain the optimal trajectory.

\begin{figure}[t]
    \centering
    \includegraphics[width=0.85\linewidth]{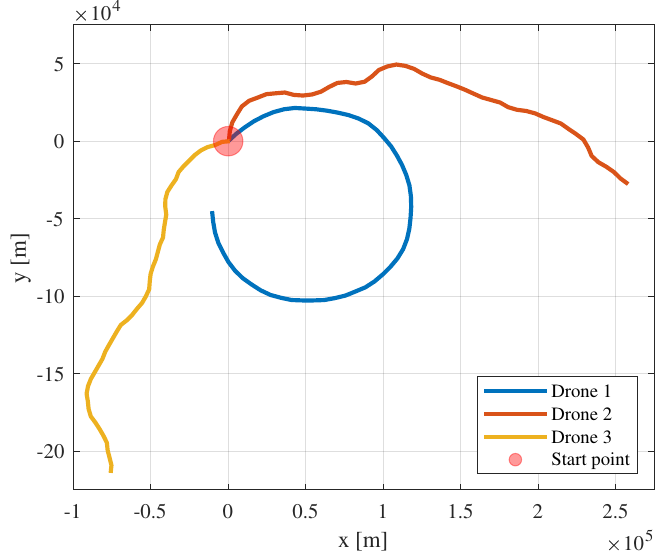}
    \caption{Tracker drone trajectories.}
    \label{fig:2}
\end{figure}

\section{Performance Analysis}
The performance of the position estimator and the path planner are evaluated, along with simulation results for single-hop and multi-hop communication. The scenario assumes a group of four drones, with one acting as a relay and the other three as trackers. Three sets of target models are available, categorized based on their speeds as slow, normal, and fast. Velocities for each drone are selected from a set of ${25,30,35}$ [m/s]. Initial positions are assigned as follows: $\mathbf{pos_{1}} = [0,0]$ [m], $\mathbf{pos_{2}} = [500,500]$ [m], and $\mathbf{pos_{3}} = [1000,0]$ [m]. 

The trajectories of the drones that are estimated are depicted in Figure \ref{fig:2}. The velocity and angle of the tracker drones change approximately every 5 minutes during the simulation. The positions of the drones are represented with data points obtained every 30 minutes.

\subsection{Estimation of Tracker Drones}
A Monte Carlo simulation is conducted to check the robustness of the estimator. The simulation time was set to 171 minutes and the following values were used for the Kalman Filter: $\mathbf{\sigma_{acc}} = [0.3,0.3]$, $\mathbf{\sigma_{gps}} = [3,3]$ and $\mathbf{\sigma_{vel}} = [0.3, 0.3]$.

\begin{figure}[t]  
    \centering
       \includegraphics[width=\linewidth]{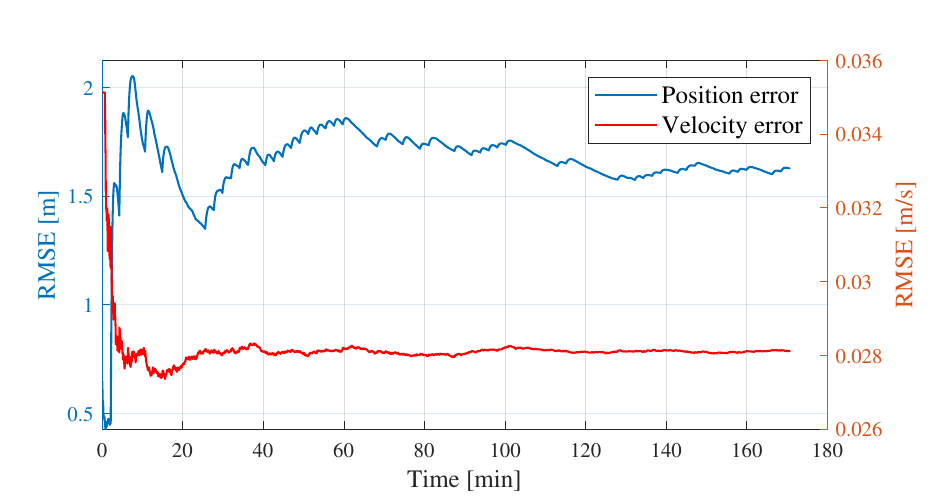}
    \caption{RMSE values for the position and velocity estimation obtained by Monte Carlo simulations.}
  \label{fig:3} 
\end{figure}

A Monte Carlo simulation is conducted with $M = 100$ to assess the robustness of the filter. The Root Mean Square Error (RMSE) values obtained for the estimated position and velocity are presented in Figure \ref{fig:3}. The plot reveals that the average RMSE values for the estimated quantities converge after some time. The error in position estimation is approximately 2 [m], whereas the error in velocity estimation is less than 0.03 [m/s].

\subsection{Path Planning}

20 different scenarios were simulated to assess the performance of the path planner for both single and multi-hop communication. 4 different communication ranges were used, namely $[50, 100, 150, 200]$ km. Each simulation lasted 171 minutes, with each time step being 2 seconds. The relay drone estimates the positions and velocities of the other drones every 30 seconds. The optimization problem is solved under two cases: with no horizon (F=1) and 8 seconds optimization horizon (F=4) window. To solve the optimization problem, constraints are imposed on the velocity and angle parameters. At each time step, the path planner selects its velocity $V_D(t) \in \{20, 30, 40\}$ [m/s] and bank angle $\theta(t) \in \{-30^{\circ}$, $0^{\circ}$, $30^{\circ}\}$. A center of mass approach is used as a baseline, where the UAV trajectory is chosen as the centroid of the three targets at time $t$. The results of the center of mass approach are treated as the baseline to compare with the optimization problem results. The calculated trajectories of the relay drone are displayed for the scenario shown in Figure \ref{fig:2}.

\begin{figure*}
    \centering
  \subfloat[Single-hop communication. \label{fig:4a}]{%
       \includegraphics[width=0.4\linewidth]{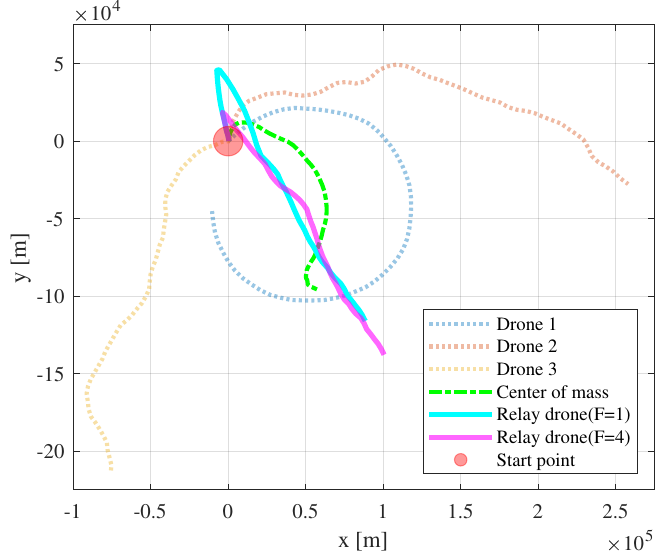}}
  \subfloat[Multi-hop communication, nearest point algorithm.\label{fig:4b}]{%
        \includegraphics[width=0.4\linewidth]{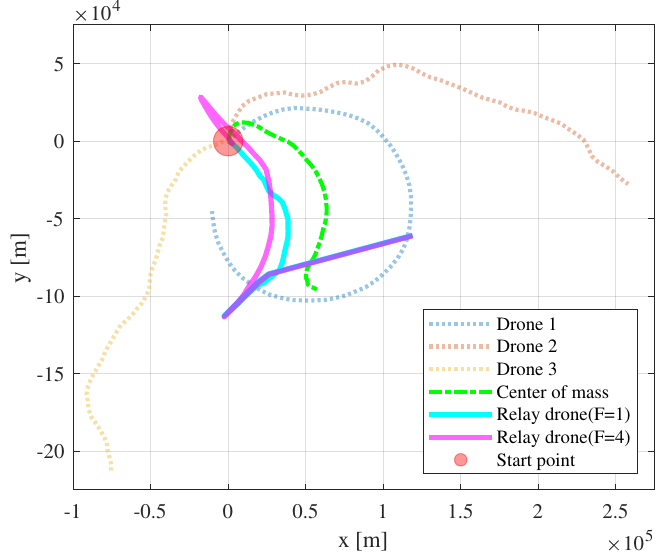}}\\
      \subfloat[Multi-hop communication, midpoint algorithm. \label{fig:4c}]{%
       \includegraphics[width=0.4\linewidth]{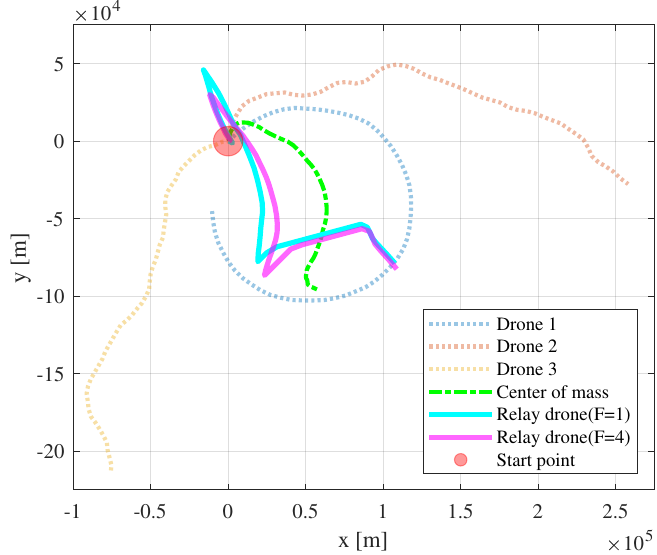}}
  \subfloat[Multi-hop communication, hybrid algorithm.\label{fig:4d}]{%
        \includegraphics[width=0.4\linewidth]{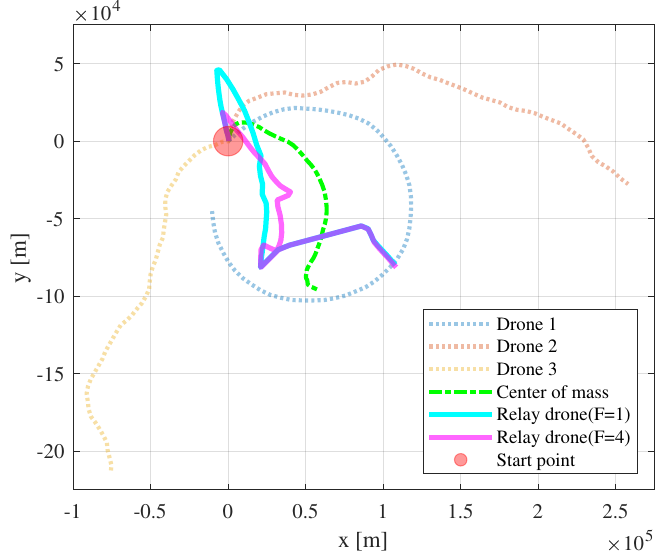}}    
    \caption{Simulation results with $R_{com} = 100$ km.}
  \label{fig:4} 
\end{figure*}

The trajectory of the relay drone optimized for the single-hop communication objective is shown in Figure \ref{fig:4a}. The movements of Drones 2 and 3 are almost symmetric on the $x=y$ line, and due to the motion of Drone 1, the center of mass trajectory has a similar shape. Conversely, the relay drone algorithm for the single-hop communication scenario attempts to establish a connection with the furthest drone. Its trajectory shows an effort to maintain a connection with Drones 2 and 3 since it is connected to Drone 1. However, since the algorithm does not take future positions into account without an optimization horizon, it fails to keep up with the tracker drones, unlike in the case where an optimization horizon is present.

The trajectory resulting from the nearest point algorithm in the multi-hop communication scenario is depicted in Figure \ref{fig:4b}. Among the drones, Drone 3 is the furthest, while Drone 1 and Drone 2 are in proximity to each other. Within this set, Drone 1 is closer to the relay drone. Thus, the relay drone stays in the centroid of Drone 1 and Drone 3, using Drone 1 to establish communication with Drone 2. As Drone 3 moves farther away, the relay drone adjusts its trajectory to the centroid of Drone 1 and Drone 2 to remain connected with both drones. However, the relay drone relies on Drone 1 to communicate with Drone 3 and does not consider the distance of Drone 3.

The trajectory of the relay drone is influenced by the midpoint algorithm, which considers the positions of all tracker drones. The optimal relay drone trajectories determined by the algorithm are depicted in Figure \ref{fig:4c}. At the start of the simulation, the trajectory is similar to the nearest point algorithm, aiming to be in the center of Drone 3, which is the furthest drone, and the centroid of Drone 1 and Drone 2. However, before Drone 2 becomes the furthest drone, the relay drone adjusts its movement in response to its movement due to the objective function. This is the main difference between the nearest point and midpoint algorithms. Towards the end of the simulation, Drone 3 becomes the furthest drone again, and the relay drone moves towards the $x=y$ line. The impact of the optimization horizon is similar to that of the single-hop communication scenario, with more intricate movements when the optimization horizon is taken into account.

In the hybrid algorithm, the trajectory at the beginning is identical to that of the single-hop communication scenario since the relay drone tries to establish a direct connection with all tracker drones. After some time, when the single-hop communication is no longer feasible, it adopts the midpoint objective, and its trajectory is almost identical to the midpoint algorithm. The impact of the horizon window is similar to the earlier scenarios. The trajectory of the relay drone is shown in Figure \ref{fig:4d}. 

At the beginning of the simulation, all drones are interconnected, and the Mayer term of the optimization determines the trajectory. When the distance between the relay drone and the tracker drone reaches the communication range limits, the Bolza term of the problem dominates the objective function. The sharp turns in the trajectory are due to the Mayer term, while the small distortions on the trajectory are results of the Bolza term. The incorporation of an optimization horizon results in increased performance of the path planner. However, there are several issues associated with introducing an optimization horizon. Firstly, the objective function becomes more complex, and the dimensionality increases exponentially with an increase in the horizon window. Moreover, the optimization operates on the estimator results, therefore the optimization horizon cannot exceed the estimation window. Longer horizon windows decrease the precision to dynamic changes in the environment.

\begin{figure}[t]
    \centering
    \includegraphics[width=\linewidth]{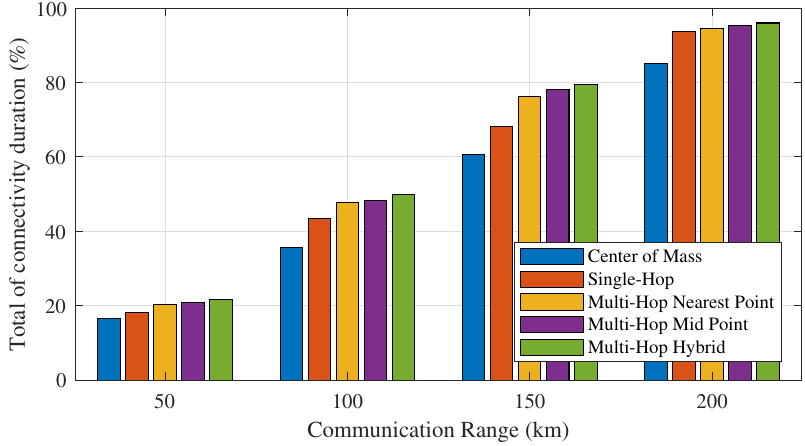}
    \caption{Total connectivity duration of the proposed algorithms.}
    \label{fig:5}
\end{figure}
Table \ref{all_results} displays the achieved connectivity duration (minutes) of the proposed algorithms in the simulation. The single-hop communication method improves the baseline technique by 7.89\% in the absence of an optimization horizon and 13.15\% with its application. The nearest point algorithm achieves a 36.84\% improvement, irrespective of the optimization horizon. The midpoint algorithm leads to an improvement of 42.10\% and 35.52\% with and without optimization horizon, respectively. Finally, the hybrid algorithm enhances the baseline by 42.10\%.

The duration of network connectivity during simulations was calculated for four different communication ranges. Figure \ref{fig:5} illustrates the correlation between the communication range of the UAVs and the connectivity time. It also demonstrates the improvement achieved by the proposed methods compared to the baseline method. For instance, when the communication range is $R_{com} =$ 100 km, the single-hop communication path planner achieves connectivity for 43.99\% of the total simulation time, while the baseline method achieves only 35.71\%. Furthermore, the multi-hop communication path planners exhibit better performance. The nearest point algorithm has a connectivity duration of 47.56\%, the midpoint algorithm has 48.73\%, and the hybrid algorithm has 50.26\%. 

The present outcomes are inadequate as key performance indicators for evaluating the planners, given that their accuracy is significantly influenced by the extent of communication coverage. In order to mitigate the impact of communication range, the duration of sustained network connectivity is calculated and subsequently compared to the period during which network connectivity is ensured by the relay drone.

\begin{table}[t]
    \renewcommand{\arraystretch}{1.1}
    \caption{Performance Evaluation of Proposed Networks on Test Set}
    \label{all_results}
    \centering
    \begin{tabular}{l c c c c}
    \hline
        \bf{Algorithm} & \bf{No Horizon} & \bf{with Horizon (8 seconds)} \\ \hline
        Single-hop & 82 & 86 \\ 
        Nearest point & 104 & 104 \\ 
        Midpoint & 103 & 108 \\ 
        Hybrid  & 108 & 108 \\ \hline
        Center of mass  & 76 & 76 \\ \hline
    \end{tabular}
\end{table}

\section{Theoretical Limits on Network Lifetime}
In this section, we investigate whether the relay UAV can be positioned such that a connected network topology is obtained among all UAVs if the positions of tracker drones are exactly known at time $t$. For each communication type, a seperate geometric problem is formulated.

\begin{figure}[t]
    \centering
    \includegraphics[width=0.8\linewidth]{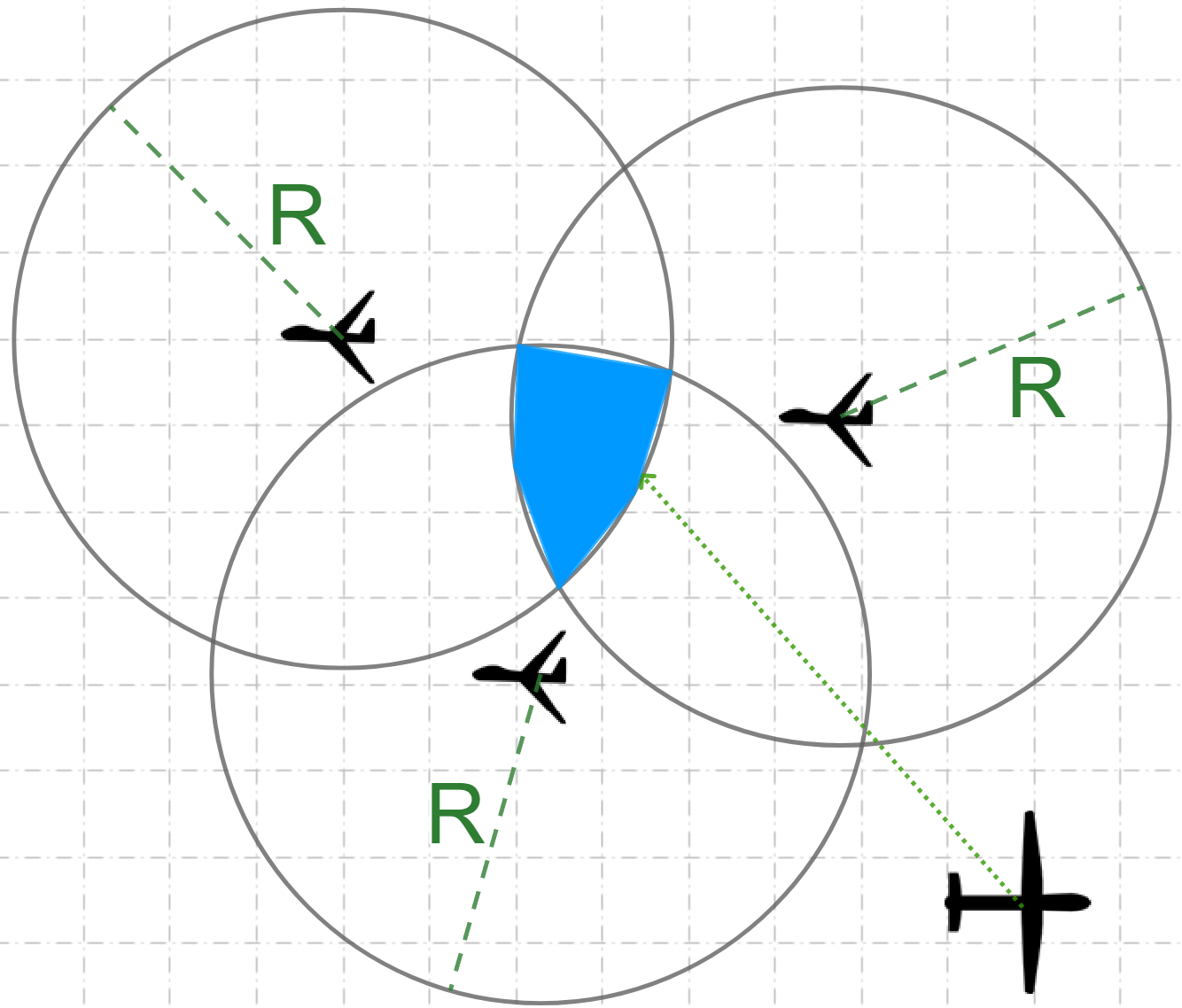}
    \caption{Connectivity area (blue) for relay drone to form a single-Hop communication.}
    \label{fig:8}
\end{figure}

\subsection{Single-Hop Communication}

To have a connected single-hop network, all tracker drones have to communicate with relay drone directly. At time t, let $d_i^{(x,y)}(t)$ denote the distance between point $(x,y)$ and the tracker drone $i$. If the distance $d_i^{(x,y)}(t) \geq R_{com}$, the connection link cannot be formed. To call a topology feasible, there at least needs to be one point $(x,y)$ in the map that satisfies $d_i^{(x,y)}(t) \leq R_{com}$  for $\forall i = 1,...,N$. 

At each time instance $t$, circles are formed centered at the positions of tracker drones with radius $R_{com}$, and the intersection area between N circles is checked. The sum of these intersection areas over the simulation time gives the maximum duration for which the network can remain connected. Figure \ref{fig:8} illustrates a feasible area for $N = 3$.

\subsection{Multi-Hop Communication}

In the case of multi-hop communication, it is necessary to have a connected network graph among the drones. However, determining the feasibility of placing the relay drone to achieve a connected network topology is more complex than in the single-hop communication case. For instance, in some cases, the relay drone must bridge two subsets of tracker drones that are already communicating within the subsets. The objective is to determine if there is a viable area for the relay drone placement in a multi-hop network with $N = 3$ tracker drones. To cover all possible cases, different topology configurations are examined. The aim is to find a feasible area that guarantees a connected topology between the drones.

\begin{figure}[t]
    \centering
    \includegraphics[width=\linewidth]{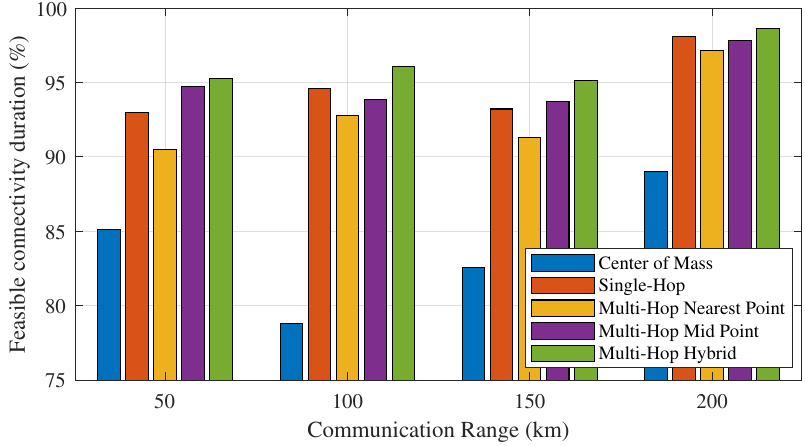}
    \caption{Relative connectivity duration of the proposed algorithms.}
    \label{fig:6}
\end{figure}

Circular areas with radius $R_{com}$ and centers $\mathbf{pos_i(t)}$ are drawn to represent the communicable areas for each drone. It is assumed that any point within the circular area can communicate with drone $i$. The distance between tracker drone $i$ and $j$ is denoted by $d_{ij}$. It is assumed, without loss of generality, that $d_{12} \leq d_{13} \leq d_{23}$. The distances between the drones are compared to the communication range to determine the existence of a feasible area The feasible points must lie within the union of the communicable areas, as the relay drone will be unable to communicate with any drones otherwise.

Three possible scenarios arise when dealing with 3 drones: all three are able to communicate, only two are able to communicate, or none are able to communicate. In the case where all three drones can communicate, any point within the union of the communicable areas of each drone is a feasible point for the relay drone, given that $d_{13} \leq R_{com}$. In this scenario, $d_{23}$ does not impact communication since the three drones are connected even without the relay drone.  An illustration of the configuration of this case can be seen in Figure \ref{fig:8}.

When only 2 drones can communicate, it implies that $d_{12} \leq R_{com}$ and $d_{13} \geq R_{com}$. This case can be further divided into two subcases based on $d_{13}$ being either between $R_{com}$ and $2R_{com}$ or larger than $2R_{com}$. In the former case, the relay drone can be placed at the intersection area of drone $1$ and drone $3$ to connect both of them, since drone $1$ and drone $2$ are already connected. Consequently, a connected network topology is formed. However, in the latter case, it is impossible to place the relay drone in a way that connects drone $1$ and drone $3$, and thus a connected network topology cannot be established.

In the event that there is no connection between the drones, where the distance between two drones is greater than or equal to the communication range, the relay drone must establish communication with all three drones to create a connected network topology. This is equivalent to the process of determining whether a single-hop network can be established between the drones, as discussed previously. For two circles to have an intersection area, the distance between their centers should be smaller than the summation of their radii. To have an intersection area between 3 circles, the distance between the circles should satisfy the following inequality $R_{com} \leq d_{12} \leq d_{13} \leq d_{23} \leq 2R_{com}$. There are two possible cases for this configuration with intersection area and no intersection area. 

A feasibility check is carried out for each $t$, and the duration that the topology can remain connected is determined by summing the results throughout the simulation time to obtain the maximum possible duration. We evaluated the maximum amount of time that the network can remain connected for each target mobility case and communication range. For each algorithm, the ratio between the time the relay drone maintains network connectivity and the maximum possible duration of network connectivity, called relative connectivity duration, is calculated.

\begin{figure}[t]
    \centering
    \includegraphics[width=\linewidth]{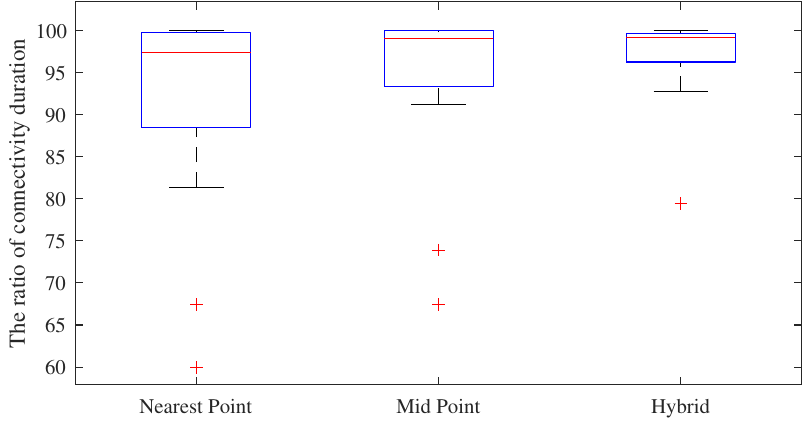}
    \caption{Box plot of multi-hop algorithms.}
    \label{fig:7}
\end{figure}

The results in Figure \ref{fig:6} show that both in single-hop and multi-hop communication cases, relative connectivity durations increase with the proposed algorithms compared with the center of mass baseline method. The algorithm for single-hop communication has better results than the nearest point algorithm since the feasible area increases in the multi-hop communication scenario. Although the nearest point algorithm planner has longer connection duration, it is not as relatively successful as the single-hop communication planner in terms of being in the connectivity area.

All three algorithms for the multi-hop communication scenario improves the performance in terms of relative connectivity duration compared to the baseline algorithm. Figure \ref{fig:6} shows that the midpoint algorithm performs better than the nearest point algorithm, and hybrid algorithm performs better than both. Figure \ref{fig:7} presents the statistics of the relative connectivity duration for all simulations for the multi-hop communication. It is observed that among the three algorithms considered, the hybrid algorithm not only achieves the longest connectivity duration, but also has the most robust performance.

\section{Conclusion}
In this paper, motion planning for a relay UAV for establishing a connected network topology is investigated, where each UAV in the network has the ability to communicate with other drones. The communication is studied under two scenarios: single-hop and multi-hop. Single-hop communication involves a direct link from the tracking UAVs to the relay UAV, which requires reliable and low delay communication to ensure effective surveillance. On the other hand, in multi-hop communication, the relay UAV can communicate with the other UAVs through intermediate nodes, which increases the duration that the UAVs remain connected to achieve higher resilience. The optimization problems are formulated for both single-hop and multi-hop communication, aiming to maximize the connected time of the network topology. 

Center of mass method which selects the trajectory of the relay UAV as the centroid of the three targets at all times is used as the baseline for comparison. The performance of the algorithm used in the single-hop communication shows that the solution with no time horizon shows 13\% improvement, while the solution with a time horizon shows 16\% improvement over the baseline. The formulation with optimization horizon provides the best results, which is expected since the trajectory is optimized with more knowledge about the tracker drones. For the multi-hop communication scenario, trajectories generated by the proposed algorithms keep the relay UAV connected for over 90\% of the maximum possible duration of having a connected topology. Amongst the proposed three algorithms, the hybrid algorithm is found to be the most effective algorithm with over 95\% of the total network lifetime achieved.  

\ifCLASSOPTIONcaptionsoff
  \newpage
\fi



\bibliographystyle{IEEEtran}
%
%
\bibliography{IEEEabrv,references.bib}
\end{document}